%% file: sample.tex
\documentclass[twoside,11pt]{article}

\usepackage{jmlr2e}
\usepackage[scr=boondoxo,scrscaled=1.05]{mathalfa}
\usepackage{amsmath}
\usepackage{amssymb}
\usepackage{thm-restate}
\usepackage{hyperref}
\usepackage{cleveref}
\usepackage{mathtools}
\usepackage{sidecap}
\usepackage[font={small}]{caption}
\jmlrheading{}{}{}{}{}{}{Shit and Ezhov et al.}
\ShortHeadings{Semi-Implicit Neural PDE Solver}{Shit and Ezhov et al.}
\firstpageno{1}

\begin{document}

\title{Semi-Implicit Neural Solver for Time-dependent \\
Partial Differential Equations}

\author{%
\name Suprosanna Shit\thanks{Contributed Equally} \thanks{Corresponding Author}  \email suprosanna.shit@tum.de\\
\addr Department of Informatics, Technical University Munich
\AND
\name Ivan Ezhov$^\ast$ \email ivan.ezhov@tum.de\\
\addr Department of Informatics, Technical University Munich
\AND
\name Leon M\"achler \email  leon-philipp.machler@ens.fr\\
\addr Département d’informatique, École normale supérieure, Paris
\AND
\name Abinav R. \email abinav@deepc.ai\\
\addr deepc GmbH, Munich
\AND
\name Jana Lipkova \email jlipkova@bwh.harvard.edu\\
\addr Brigham and Women's Hospital, Harvard Medical School
\AND
\name Johannes C. Paetzold   \email johannes.paetzold@tum.de\\
\addr Department of Informatics, Technical University Munich
\AND
\name Florian Kofler   \email florian.kofler@tum.de\\
\addr Department of Informatics, Technical University Munich
\AND
\name Marie Piraud \email marie.piraud@helmholtz-muenchen.de\\
\addr Helmholtz AI, Helmholtz Zentrum München
\AND
\name Bjoern H. Menze \email bjoern.menze@uzh.ch\\
\addr Department of Quantitative Biomedicine, University of Zurich}

% \editor{Kevin Murphy and Bernhard Sch{\"o}lkopf}
\editor{}
\maketitle

\begin{abstract}%   <- trailing '%' for backward compatibility of .sty file
Fast and accurate solutions of time-dependent partial differential equations (PDEs) are of pivotal interest to many research fields, including physics, engineering, and biology. Generally, implicit/semi-implicit schemes are preferred over explicit ones to improve stability and correctness. However, existing semi-implicit methods are usually iterative and employ a general-purpose solver, which may be sub-optimal for a specific class of PDEs. In this paper, we propose a neural solver to learn an optimal iterative scheme in a data-driven fashion for any class of PDEs. Specifically, we modify a single iteration of a semi-implicit solver using a deep neural network. We provide theoretical guarantees for the correctness and convergence of neural solvers analogous to conventional iterative solvers. In addition to the commonly used Dirichlet boundary condition, we adopt a diffuse domain approach to incorporate a diverse type of boundary conditions, e.g., Neumann. We show that the proposed neural solver can go beyond linear PDEs and applies to a class of non-linear PDEs, where the non-linear component is non-stiff. We demonstrate the efficacy of our method on 2D and 3D scenarios. To this end, we show how our model generalizes to parameter settings, which are different from training; and achieves faster convergence than semi-implicit schemes.
\end{abstract}


\input{01_intro}
\input{02_methodology}
\input{03_experiments}

\section{Conclusion}
This work introduces a novel implicit neural scheme to solve time-dependent PDEs in arbitrary geometry and boundary conditions. We leverage an existing semi-implicit update rule to design a learnable iterator that provides theoretical guarantees. The learned iterator achieves faster convergence compared to the existing semi-implicit solver and produces a more accurate solution for a fixed computation budget. Importantly, we empirically demonstrate that training on a single parameter setting is enough to generalize over other parameter settings which confirms our theoretical results. The learner neural solver offers computationally and scalable alternative to standard numerical approaches. The increased computational efficiency expand computational possibilities  and enables simulations of complex real world systems.

\section*{Acknowledgement}
S. Shit and I. Ezhov  are  supported  by  the  Translational  Brain  Imaging Training Network under the EU Marie Sklodowska-Curie programme (GrantID:  765148). We thank Prof. Dr. Elisabeth Ullmann for giving critical feedback on this work.

\clearpage

\appendix
\input{04_appendices}

\clearpage
\bibliography{sample}

\end{document}

%% file: 01_intro.tex
\section{Introduction}
\label{sec:intro}
Time-dependent partial differential equations (PDEs) are an essential mathematical tool to describe numerous physical processes in various disciplines, such as wave propagation \citep{zhou2004periodic}, quantum transport \citep{manzano2012quantum,ezhov2016influence}, cell diffusion \citep{hinderliter2010isdd,lipkova2019personalized}, among others. Solving the initial-value problem and the boundary-value problem accurately in a computationally efficient way is the primary research interest for these PDE problems. 

The numerical solution of time-dependant PDEs relies on the appropriate spatio-temporal discretization. Spatial discretization can be implemented using a finite difference, finite element, or finite volume method. For temporal discretization under Eularian settings, either explicit, implicit, or semi-implicit methods can be used. Explicit temporal update rules are generally a single or few forward computation steps. In contrast, implicit or semi-implicit update rules, such as Crank-Nicolson's scheme, require a fixed-point iterative solver. In all the methods mentioned above, smaller time steps and finer spatial resolution facilitate a more accurate solution. At the same time, it substantially increases the computational burden. Moreover, the maximum allowed spatio-temporal resolution is also upper bounded by numerical stability criteria. It can be observed that, in contrast to explicit methods, implicit and semi-implicit methods offer relaxed stability constraints (sometimes unconditionally stable) for admissible time steps at the expense of an increased computational cost caused by the iterative solver.

In recent times, the use of neural networks, e.g., as proposed by \cite{raissi2019physics} has gained significant attention for supporting numerical computations. Superior performance is achieved in solving forward simulations \citep{magill2018neural,li2020fourier, tompson2017accelerating,ezhov2020real,kochkov2021machine,greenfeld2019learning,shit2021velocity} and inverse-problems \citep{ref_article13,ref_article12,long2017pde,ezhov2019neural,greenberg2019automatic}. On the contrary to the well understood and theoretically grounded classical methods, the deep learning-based approaches rely mainly on empirical validity. Recently, \cite{hsieh2018learning} developed a promising method to learn numerical solvers while providing a theoretical convergence guarantee. They demonstrate that a feed-forward fully-convolutional network (FCN) trained to correct the error of a single iteration of a linear solver can deliver a faster solution than the hand-crafted solver. Astonishingly, for time-dependent PDEs, the temporal update step of all previously mentioned applications neural schemes relies on an explicit forward Euler method; none of them is exploiting the powerful implicit and semi-implicit methods. Incorporating implicit iterative schemes into the neural solvers would potentially expand the applicability of neural architectures in real-world applications by solving time-dependent PDEs much more efficiently.

\subsection{Our Contribution}
In this paper, we introduce the first neural solver for time-dependant PDEs. First, we construct a neural iterator from a semi-implicit update rule for linear PDEs with the Dirichlet boundary condition. Then, we replace a single iteration of the semi-implicit scheme with a learnable parameterized function such that the fixed point of the algorithm is preserved. We provide theoretical guarantees, similar to \cite{hsieh2018learning}, which prove that our proposed solution converges to the correct solution and generalizes over parameter settings very different from the ones seen at training time. Subsequently, we show how the method can be extended to other types of BC through diffuse domain approximation \citep{li2009solving}, which we illustrate on Neumann boundary condition. Next, we extend our learnable framework beyond linear PDEs to a class of non-linear PDEs. We validate our method on 2D and 3D experiments with Dirichlet, Neumann boundary conditions, and non-linear PDEs. Empirically, we show that our model also generalizes well across different complex geometries and produces a more accurate solution than the existing semi-implicit solver while taking much lesser computational cost.

%% file: 02_methodology.tex
\section{Method}
\label{sec:method}

In the following, we present the iterative semi-implicit scheme (Section 2.1), our proposed learnable approximation of the iterative update (Section 2.2), how to implement the Neumann boundary condition (Section 2.3), demonstrate its applicability to a broad class of non-linear PDEs (Section 2.4), and describes the neural solver learning setup (Section 2.5).

\subsection{Time-dependent Linear PDEs with Dirichlet Boundary Condition}

First, we consider the initial value problem of variable of interest $\mathscr{u}$ governed by a time-dependent linear PDE with Dirichlet boundary condition of the following form,
\begin{equation}
	\frac{\partial \mathscr{u}}{\partial t} = \mathcal{F}(\mathscr{u};\{\partial_i, \mathscr{a}_i \}_{i=1:N}), \forall x \in \Omega, \text{ s.t. } \mathscr{u}(x,t) = \mathscr{b}(x,t), \forall x \in \Gamma \text{ and } \mathscr{u}(x,t_0) = \mathscr{u}_0 \label{eq:linear_dirichlet}
\end{equation} 

\noindent where $\Omega$ is the domain (e.g., $\in \mathbb{R}^2$ or $\in \mathbb{R}^3$), $\Gamma$ is the boundary with boundary values $\mathscr{b}(x,t)$ and $\mathscr{u}_0$ the initial value at time $t=t_0$. $\mathcal{F}(\mathscr{u},\{\partial_i, \mathscr{a}_i \}_{i=1:N})$ is a linear combination of spatial partial differential operators $\{ \partial_i \}_{i=1:N}$ and its corresponding parameter $\{ \mathscr{a}_i \}_{i=1:N}$. Without the loss of generalizability, we choose a uniform discretization step $\Delta x$ for all spatial dimensions and we can discretize $\mathcal{F}$ it in the following matrix form: 
\begin{equation} 
    F(u;\mathbb{D},\mathbb{A}) = \sum_{i = 1}^{N} \frac{A_iD_i}{\Delta x^{p_i}}u
\end{equation}
where $\mathbb{D}=\{ D_i \}_{i=1:N}$, $\mathbb{A}=\{ A_i \}_{i=1:N}$, $A_i$ is a diagonal matrix consisting of the values of $\mathscr{a}_i$ corresponding to the discrete differential operator $D_i$ of order $p_i$, which is a Toeplitz matrix. We denote $u$ at time $t$, as $u_t$. A first order semi-implicit update rule to get $u_{t+\Delta t}$ from $u_{t}$ (with time step $\Delta t$) is given by
\begin{equation}
	\frac{u_{t+\Delta t} - u_t}{\Delta t} = \theta F(u_{t+\Delta t};\mathbb{D},\mathbb{A}) + (1- \theta)F(u_{t};\mathbb{D},\mathbb{A}) \quad [0 < \theta \leq 1]
\end{equation}
To obtain $u_{t+\Delta t}$, one needs to solve the following linear system of equations
\begin{equation}
	\left( I - \theta\Delta t \sum_{i = 1}^{N}\frac{A_i d_i}{\Delta x^{p_i}} \right) u_{t + \Delta t} = \theta\Delta t \sum_{i = 1}^{N}\frac{A_i (D_i - d_i I)}{\Delta x^{p_i}} u_{t + \Delta t} + c (u_t, \mathbb{D}, \mathbb{A}, \Delta x, \Delta t)\label{eq:implicit_update}
\end{equation}
where $c$ is independent of $u_{t + \Delta t}$ and $d_i$ is the central element of the central difference discretization of $D_i$. Note that for central difference scheme, $D_i-d_i I$ is real, zero-diagonal, and either circulant or skew-circulant matrix. 

One can use an iterative scheme to compute $u_{t+\Delta t}$ from an arbitrary initialization $u^0$ on the right-hand-side of Eq. \ref{eq:implicit_update}. We denote $m^{th}$ updated value of $u_{t+\Delta t}$ as $u^m$, and for the ease of notation, we introduce $\{\Lambda_i\}_{i=1:N}$ and $\tilde{c}$ as following:
\begin{align}
	\Lambda_i &= \left( I - \theta\Delta t \sum_{j = 1}^{N}\frac{A_j d_j}{\Delta x^{p_j}} \right)^{-1} \frac{\theta\Delta t A_i}{\Delta x^{p_i}}\\
	\tilde{c} &= \left( I - \theta\Delta t \sum_{j = 1}^{N}\frac{A_j d_j}{\Delta x^{p_j}} \right)^{-1} \left( I + (1 - \theta)\Delta t \sum_{i = 1}^{N} \frac{A_i D_i}{\Delta x^{p_i}} \right) u_t
\end{align}
Using the $\{\Lambda_i\}_{i=1:N}$ and $\tilde{c}$ notations the iterator can be written as:
\begin{equation}
    	u^{m+1} =  \sum_{i = 1}^{N} \Lambda_i (D_i - d_i I)u^m + \tilde{c},
\end{equation}
and by enforcing the Dirichlet boundary condition using a projection step with a binary boundary mask $G$, the iterator becomes:
\begin{equation}
    	u^{m+1} = G \left( \sum_{i = 1}^{N} \Lambda_i (D_i - d_i I)u^m + \tilde{c} \right) + (I - G)b_{t+\Delta t} \label{eq:new_it}
\end{equation}
where $b_{t+\Delta t}$ is the boundary value. Notice that $\tilde{c}$ is independent of $u^m$. In \cite{hsieh2018learning}, a formal definition of a \textit{valid iterator} for solving linear system of equation is provided invoking convergence and fixed point criteria. Using the same definition, we show that with an appropriate choice of $\Delta x, \Delta t$ and $\theta$, Eq. \ref{eq:new_it} is a \textit{valid iterator}.

\begin{restatable}[]{theorem}{thmone}
\label{thm1}
	For an appropriate choice of $\Delta x, \Delta t$ and $\theta$, the linear iterator
	
\begin{equation}
u^{m+1} = G \left( \sum_{i = 1}^{N} \Lambda_i (D_i - d_i I)u^m + \tilde{c} \right) + (I - G)b_{t+\Delta t}\nonumber
\end{equation} 
is valid.
\end{restatable}
\begin{proof}
See Appendix.
\end{proof}
For simplicity of notation, we drop the superscript of $u$ for a single update. Hence right hand side of the Eq. \ref{eq:new_it} can be seen as a linear operator 
\begin{equation}
\Psi(u)= Lu + k \label{eq:iterator}
\end{equation}

where $L=G \left( \sum_{i = 1}^{N} \Lambda_i (D_i - d_i I)\right)$ and $k=G\tilde{c}+(I-G)b_{t+\Delta t}$. 

\subsection{Neural Solver}
We propose the following end-to-end trainable iterator to replace the analytically derived iterator in Eq. \ref{eq:iterator}, following a similar structure to \cite{hsieh2018learning}
\begin{equation}
	\Phi_H(u) = \Psi(u) + G \left( \sum_{i=1}^{N} \Lambda_i H_i w \right) \label{eq:neural_iter}
\end{equation}
where $w = \Psi(u) - u$ and $H_i$ are the learnable operators, which satisfy $H_i 0 = 0, \forall i \leq N$. Substituting $w$ in Eq. \ref{eq:neural_iter}, we get to the linear neural operator
\begin{equation}
    \Phi_{H}(u) = L' u+k'
\end{equation}
where $k'$ is independent of $u$, and $L'=L+G\sum_{i=1}^{N}\Lambda_i  H_i(L-I)$. Following, we show that the neural iterator has the same fixed point ($u^*: \Phi_{H}(u^*)=u^*$) as the hand-designed iterator in Eq. \ref{eq:new_it}.

\begin{restatable}[]{lemma}{lemmatwo}
\label{lemma2}
	For a given linear PDE problem $( \mathbb{D},\mathbb{A}, G, u_t, b_{t+\Delta t}, \Delta x, \Delta t, \theta)$ and any choice of $\mathbb{H}=\{H_i\}_{i \leq N}$ that satisfies $H_i0 = 0, \forall i=1:N$, a fixed point $u_{t+\Delta t}^*$ of $\Psi$ is also a fixed point of $\Phi_H$.
\end{restatable}

\begin{proof}
Follows trivially from the definition.
\end{proof}

\noindent By construction, the neural iterator in Eq. \ref{eq:neural_iter} inherits the properties of the iterator proposed by \cite{hsieh2018learning}. The most notable one is that, if $H_i=0,~\forall i=1:N$ then $\Phi_{H}=\Psi$. Furthermore,
if $H_i=(D_i - d_i I),~\forall i=1:N$, then since $G L=L$
\begin{equation}
    \Phi_{H}(u)=\Psi(u)+G L(\Psi(u)-u)=L \Psi(u)+k=\Psi^{2}(u)
\end{equation}
which is equal to two iterations of $\Psi$. Since computing $\Phi$ requires two separate convolutions: i) $L$, and ii) $\mathbb{H}$. Since, a convolution operation needs roughly the same computation time as a differential operator, one iteration of $\Phi_{H}$ requires same order complexity of two iterations of $\Psi$. This shows that we can learn a set of $\mathbb{H}$, such that the proposed iterator $\Phi_{H}$ performs at least as good as the standard solver $\Psi$.

Following the theoretical framework of \cite{hsieh2018learning}, we extend their results for the case of the proposed neural iterator. The usefulness of any trainable iterator naturally depends on its training efficacy and the ability to generalize. Thus for trainability, we first show that the spectral norm of $\Phi_H(u)$ is a convex function of $\mathbb{H}$ and the relevant space of $\mathbb{H}$ is a convex open set. Subsequently, we show that the iterator generalizes over arbitrary initialization, boundary value, and PDE parameters.

\begin{restatable}[]{theorem}{thmthree}
\label{thm3}
For fixed $\mathbb{D},\mathbb{A}, G, u_t, b_{t+\Delta t}, \Delta x, \Delta t,\mbox{ and } \theta$, the spectral norm of $\Phi_H$ is a convex function of $\mathbb{H}$, and the set of $\mathbb{H}$, such that the spectral norm of $\Phi_H(u) < 1$ is a convex open set.
\end{restatable}
\begin{proof}
See Appendix.
\end{proof}

\noindent To prove generalizability, we first need an upper bound on the spectral norm of $\Phi_H$.

\begin{restatable}[]{lemma}{lemmafour}
\label{lemma4}
	For a choice of $\Delta x, \Delta t$ and $\theta$, if $L$ is a valid iterator, $L'$ is also a valid iterator. 
	
\end{restatable}
\begin{proof}
See Appendix.
\end{proof}
In stark contrast with previous work by \cite{hsieh2018learning}, we have several sets of parameters $\mathbb{A}, \Delta x, \Delta t$, and $\theta$ attached to the PDEs governing equation. The following theorem allows us to train the proposed model on finite parameter settings and still generalizes well.

\begin{restatable}[]{corollary}{corofive}
\label{coro5}
	For a fixed $\mathbb{D}, G$ and $\mathbb{H}$, and some $u_t', b_{t+\Delta t}', \mathbb{A}, \Delta x', \Delta t',$ and $\theta'$, if $\Phi_H(u)$ is a valid iterator for the PDE problem $(\mathbb{D},\mathbb{A}, G, u_t', b_{t+\Delta t}', \Delta x', \Delta t', \theta')$, then for all $u_t$ and $b_{t+\Delta t}$, the iterator $\Phi_H(u)$ is a valid iterator for the PDE problem $(\mathbb{D},\mathbb{A}, G, u_t, b_{t+\Delta t},$ $ \Delta x,$ $\Delta t, \theta)$, if $\Delta x, \Delta t$ and $\theta$ are chosen such that $||\Lambda_i|| < \frac{1}{\sum_{j = 1}^{N} || (D_j - d_j I)||}, \forall i=1:N$
\end{restatable}
\begin{proof}
See Appendix.
\end{proof}

\subsection{Neumann Boundary Condition}
\label{sec:phase_field}

\cite{li2009solving} proposed diffuse domain approach to accurately approximate a Neumann boundary condition on a complex geometry into a Dirichlet boundary condition on a simple geometry. They introduce a phase field function to take care of a smooth transition between the domain and boundary. We consider the following Neumann boundary condition
\begin{equation}
\frac{\partial \mathscr{u}}{\partial t}= \mathcal{F}(\mathscr{u};\{\partial_i, \mathscr{a}_i \}_{i=1:N}), \forall x \in \Omega, \text{ s.t. }  \nabla\mathscr{u}(x,t)\cdot \eta= \mathscr{b}(x,t), \forall x \in \Gamma \text{ and } \mathscr{u}(x,t_0) = \mathscr{u}_0 \label{eq:neumann_eq}
\end{equation}

where $\eta$ is the unit outwards surface normal at the boundary.
\begin{SCfigure}[][h]
\includegraphics[width=0.5\textwidth]{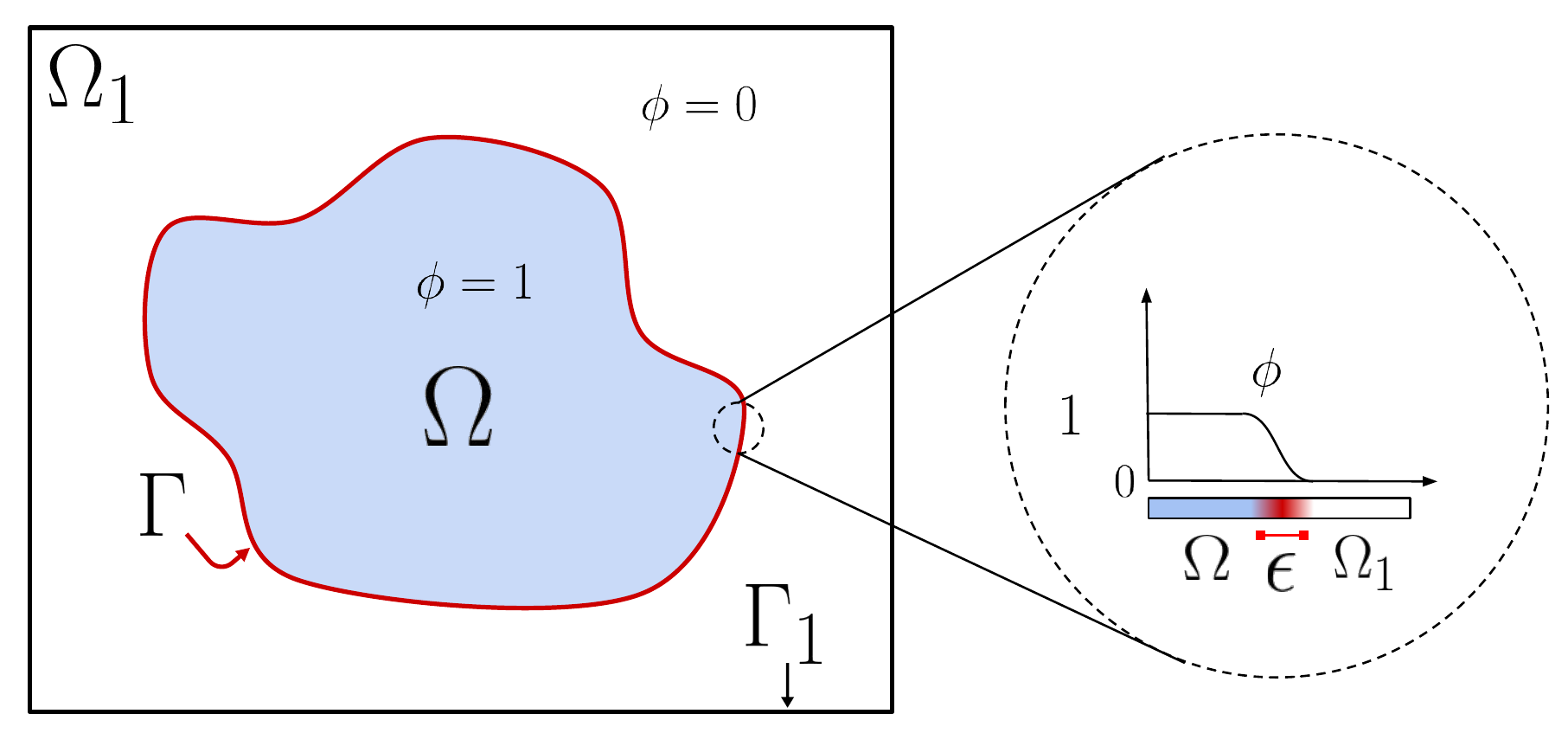}
% \vspace{-1.5em}
\caption{An illustration of the phase-field function on a complex boundary. The original BC was Neumann at the boundary $\Gamma$. Diffused domain approximation transforms it into a Dirichlet BC at the modified simpler boundary $\Gamma_1$. The smoothness factor ($\epsilon$) controls the approximation error vs robustness to noisy boundary trade-off.} \label{fig:ddm}
\end{SCfigure}

\noindent As shown in Fig \ref{fig:ddm}, the computational domain can be realized as a phase-field function $\phi$, where $\phi=1$ inside the domain and $\phi=0$ outside the domain. At the boundary, this creates a smooth transition from $1\rightarrow 0$.

Now, we solve the following approximated PDE problem 
\begin{align}
    \phi\frac{\partial \mathscr{u}}{\partial t} = \phi\mathcal{F}(\mathscr{u};\{\partial_i, \mathscr{a}_i \}_{i=1:N})+ \nabla \phi \cdot \nabla \mathscr{u}+\mathscr{b}|\nabla \phi|, \forall x \in \Omega_1\nonumber \\ \text{ s.t. } \mathscr{u}(x,t) = 0, \forall x \in \Gamma_1 \text{ and } \mathscr{u}(x,t_0) = \mathscr{u}_0
    \label{eq:ddm_eq}
\end{align}

\noindent It can be shown that without the smoothness in the boundary of phase-field function, Eq. \ref{eq:ddm_eq} reduces to Eq. \ref{eq:neumann_eq} (since within domain $\phi=1,\nabla\phi= 0$, we get the original PDE back, and at boundary $\phi=0,\nabla\phi=-\eta\neq 0$ gives the original boundary condition). The advantage of the phase-field approximation is that we can only approximate the exact surface normal orientation up to a certain accuracy for voxel-based computation, depending on the discretization resolution. Thus our surface normal orientation is noisy, resulting in an unstable boundary condition for a highly irregular shape of boundary that is a common feature of many applications. However, when we use the smoothness on the boundary, we have a consistent orientation of the resultant surface normals at the cost of the reduced magnitude of the surface normal. Eq. \ref{eq:ddm_eq} is also of linear form and can easily be converted to the solution formulation of Eq.\ref{eq:new_it} and, hence, we can apply the proposed neural solver.

\begin{figure}[t!]
	\centering
	\includegraphics[width=0.85\textwidth]{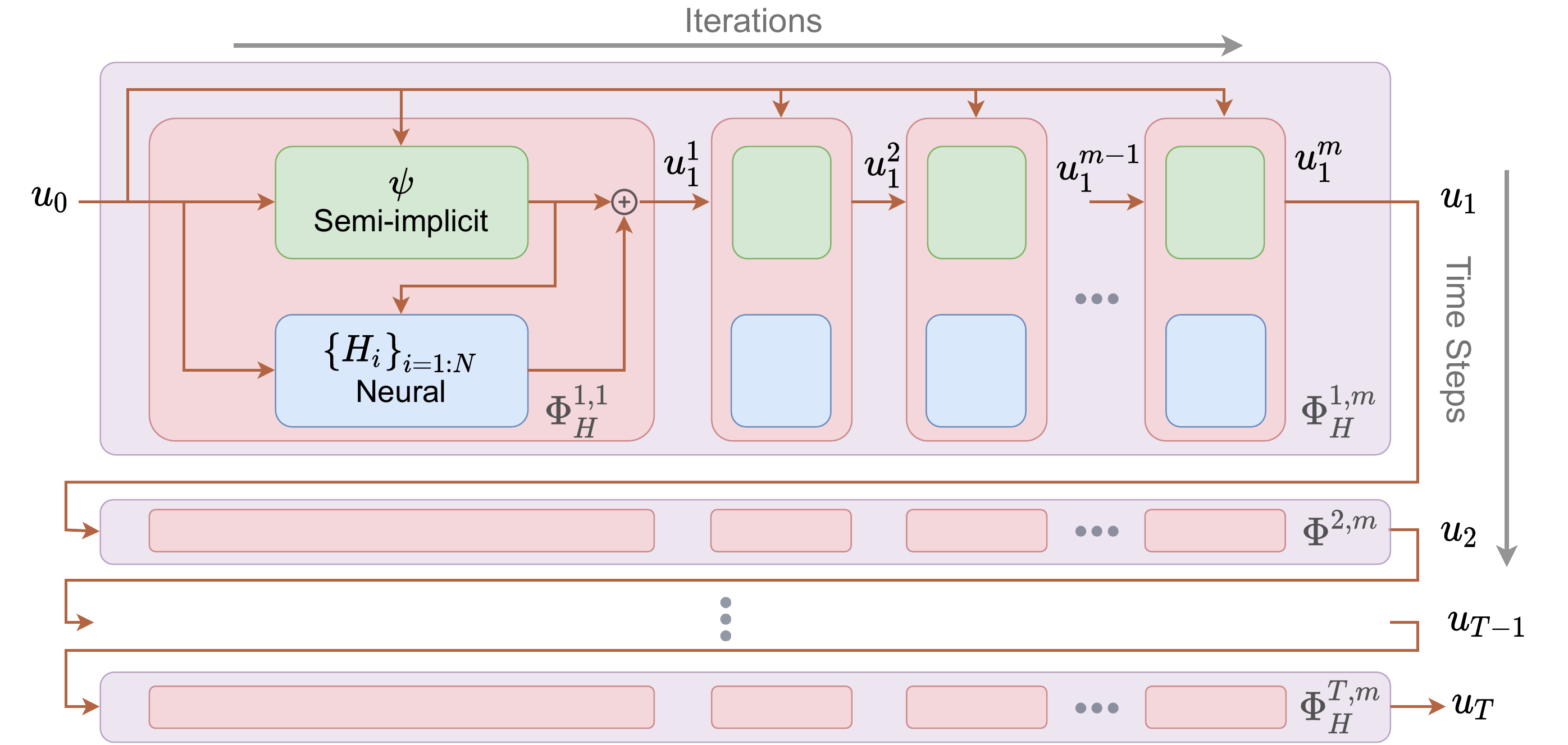}
	\caption{Semi-implicit Neural Solver: The initial condition ($u_0$) goes through a series of iterations, which consist of a semi-implicit solver ($\Psi$) and a learned neural correction ($\{H_i\}_{i=1:N}$) to produce the next solution ($u_1$). These procedure repeat at all consecutive time steps to generate the corresponding solutions ($u_1,u_2,\cdots, u_{T}$).}
	\label{fig:overview}
\end{figure}

\subsection{Extending to a Class of Non-linear PDEs}
We consider the following non-linear class of PDEs described in  \cite{boscarino2016high}:
\begin{align}
	\frac{\partial \mathscr{u}}{\partial t} = \mathcal{F}(\mathscr{u};\{\partial_i, \mathscr{a}_i \}_{i=1:N_1})+ \mathcal{G}(\mathscr{u};\{\partial_i, \mathscr{z}_i \}_{i=1:N_2}), \nonumber\\
	\forall x \in \Omega, \text{ s.t. } \mathscr{u}(x,t) = \mathscr{b}(x,t), \forall x \in \Gamma \text{ and } \mathscr{u}(x,t_0) = \mathscr{u}_0 \label{eq:non_linear}
\end{align}

\noindent where $\mathcal{F}$ is any linear operator and $\mathcal{G}$ is the non-stiff non-linear operator. Since $\mathcal{G}$ is non-stiff, we can use explicit schemes to solve the non-linear component., while we continue to use the semi-implicit formulation for $\mathcal{F}$. Thus, an implicit/explicit scheme arises [\cite{boscarino2016high}], which is used in many problems, including convection-diffusion equations, reaction-diffusion equations, collisional kinetic equations, etc. Because of using an explicit scheme for $\mathcal{G}$, it is evident that the non-linear term remains in the form of a constant in the iterative update rule. Thus, the theoretical guarantees are valid for this particular class of non-linear PDEs.

\subsection{Learning Setup}

We model each of $\{ H_i \}_{i =1: N}$ with a separate three layer convolutional neural network without any bias and non-linear activation function (c.f. Fig. \ref{fig:overview}). For each time step, the network takes the initial variable $u_0$, PDE parameters and the hyperparameters as input and produces the solution for the next time points using a fixed number of iterative forward passes. $m$ iteration for a single time step denoted as $\Phi_H(\Phi_H \dots (\Phi_H)) = \Phi^{1,m}_H$. Starting from $u_0$, we denote the solution at $t=t_0+\Delta t$ after $m$ iterations, i.e., $u_1$ as $\Phi_H^{1,m}(u_0)$. For $n$ number of time steps at $t=t_0+T\Delta t$ this process repeats and we get the solutions, i.e., $u_T$ as $\{\Phi_H^{T,m}(u_0)\}$. The reference solutions $\{u^*_{ t}\}_{t=1:T}$ can be easily obtained from $\Psi^{t,m}(u_t)$ using a sufficiently large $m$ till machine precision convergence.

We minimize the following mean squared loss function,
\begin{equation}
    \mathcal{L}=\frac{1}{T}\sum_{t=1}^{T}\left\|\Phi^{t,m}_H(u_0)-u^*_{t}\right\|; m\sim\mathcal{U}[M_1,M_2] \label{eq:costfunc}
\end{equation}
where $T$ is the maximum number of time-steps. We select the number of iteration from an uniform distribution between $[M_1, M_2]$ to avoid overfitting during training.

%% file: 03_experiments.tex
\section{Experiments}

In this section, we present experimental validation for the proposed neural solver. We aim to answer the following questions:
\begin{itemize}
    \item[Q1] How well does the neural solver converge for different sets of PDE parameters?
    \item[Q2] How good is the proposed solver for different classes of linear PDEs in different dimensionality?
    \item[Q3] Does the neural solver generalize over different computation domains, which are not seen during the training period?
    \item[Q4] How is the performance for the class of non-linear PDEs that we describe in the method section?
    \item[Q5] How is the accuracy and runtime gain for the neural solver compared to the traditional numerical schemes?
\end{itemize}

To answer these questions, we design two experiments: in 2D and 3D space, with linear and non-linear PDEs, with different PDE parametrization, in different simulation domains, and boundary conditions. In Sec. \ref{exp:2d} and \ref{exp:3d}, we detail two experiments. From this, we draw precise answers to Q1-5 and a brief discussion of the results in Sec. \ref{exp:disc}.

\subsection{2D Linear PDE with Dirichlet Boundary}
\label{exp:2d}
We consider a 2-D advection-diffusion equation of the following form
\begin{eqnarray}
    \frac{\partial \mathscr{u}}{\partial t} = [v_x, v_y] \cdot \begin{bmatrix}
\partial_{x}\mathscr{u}\\
\partial_{y}\mathscr{u}
\end{bmatrix} + [\kappa_{xx}, \kappa_{yy}] \cdot  \begin{bmatrix}
\partial_{xx}\mathscr{u}\\
\partial_{yy}\mathscr{u}
\end{bmatrix} \mbox{; subject to }\mathscr{u}(x,t_0) = \mathscr{u}_0 \label{eq:exp1}
\end{eqnarray}
where $[v_x, v_y]$ and $[\kappa_{xx}, \kappa_{yy}]$ are advection velocity and diffusivity respectively.
\begin{figure}[t!]
\centering
\includegraphics[width=0.98\textwidth]{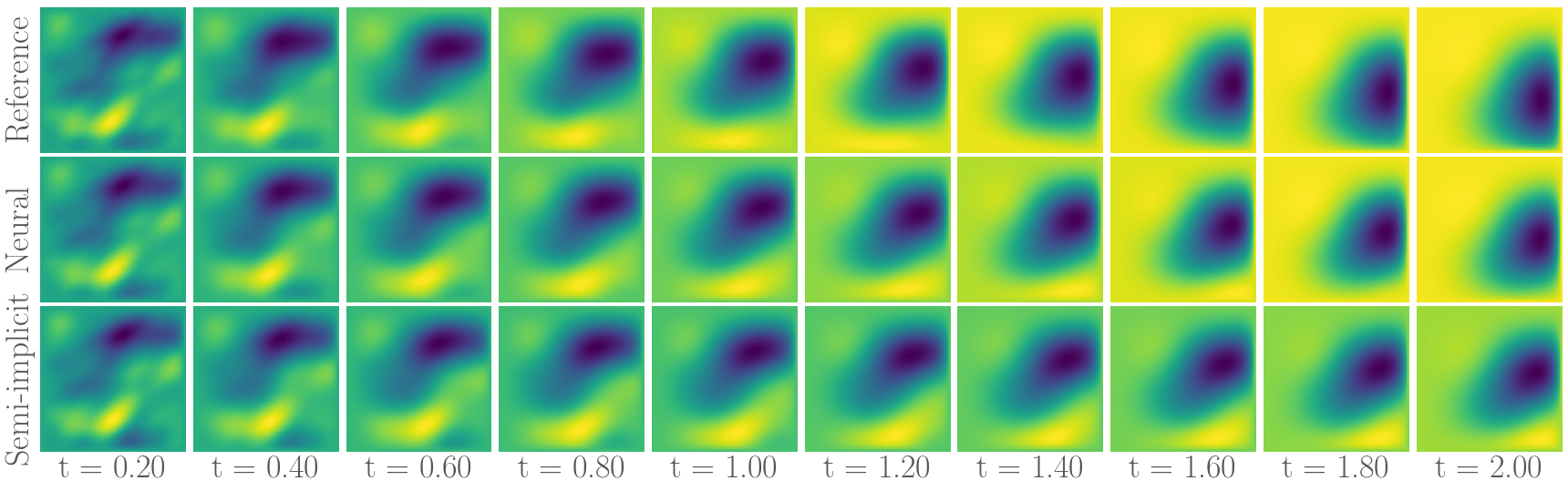}
\includegraphics[width=0.98\textwidth]{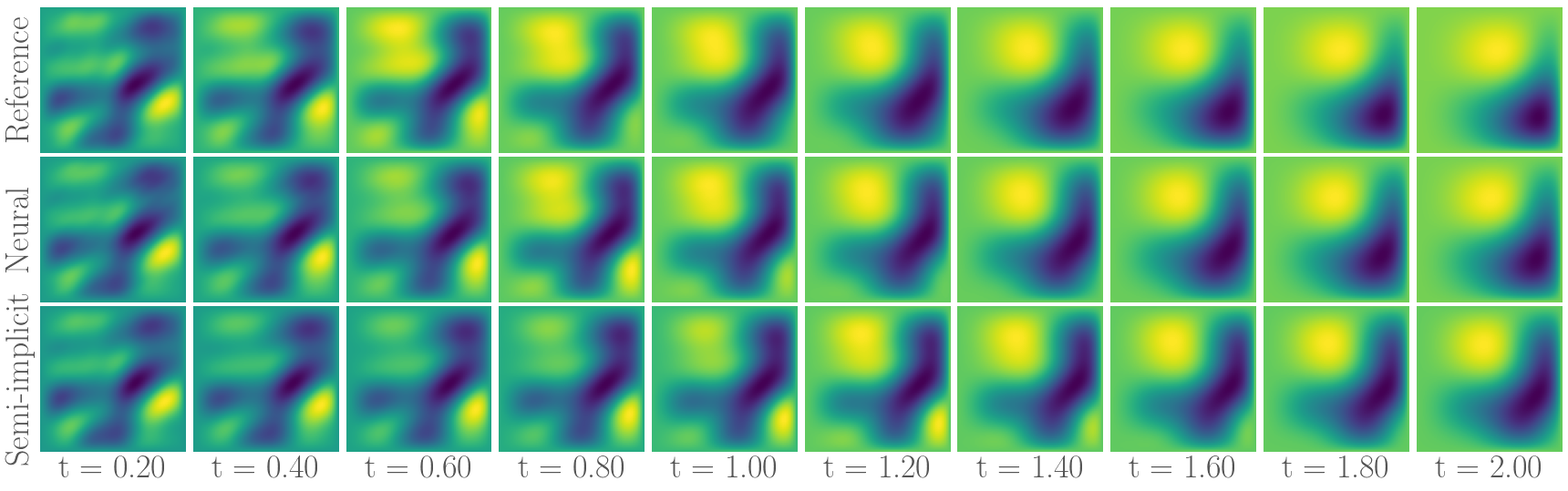}
\caption{Qualitative comparison of $u$ (c.f. Eq \ref{eq:exp1}) from the neural scheme (10 iterations) and a semi-implicit scheme (25 iterations) against the FEniCS solution for a test sequence of 10 time points. All methods use the same initial- and boundary condition. The neural update shows consistently faster convergence than semi-implicit one.}
\label{fig:fig2_exp1}
\end{figure}
\subsubsection{Data Generation}
We follow the experimental setup by \cite{long2017pde} and consider a rectangular domain of $\Omega = [0,2\pi]\times [0,2\pi]$. Elements of $[v_x, v_y]$ and $[\kappa_{xx}, \kappa_{yy}]$ are drawn from a uniform distribution of $\mathcal{U}[-2.0,2.0]$ and $\mathcal{U}[0.2,0.8]$ respectively. The computational domain is discretized by a 64 x 64 regular mesh. We assume zero Dirichlet boundary condition and the initial value is generated according to \cite{long2017pde} as $ u_0 = \lambda cos(kx+ly) + \gamma sin(kx+ly)$ where $\gamma$ and $ \lambda$  are drawn from a normal distribution of $\mathcal{N}(0,0.02)$, and, $k$ and $l$ have random values from a uniform distribution of $\mathcal{U}[1,9]$. We generate 200 simulations, each with 50 time steps, using FEniCS \citep{alnaes2015FEniCS} for $\Delta t = 0.2$. FEniCS uses finite element discretization and we found that the FEniCS solution is identical to the one we obtain from our semi-implicit solver with sufficiently large number of iterations per time step. This confirms that there is no discrepancy in the discretization and implementation of the boundary conditions for our semi-implicit and neural solver.  An exemplary time series of a test data is shown in Fig \ref{fig:fig2_exp1}.
% Since we use the implicit scheme for the time discretization, we know that it is stable for any number of time steps.

\subsubsection{Experimental Details}
We split our train, test, and validation set of the simulated time series in $80\%:10\%:10\%$. During training, we fixed the following parameters as follows $\Delta x=0.098, \Delta t=0.2, \theta=0.9$. The elements of $[v_x, v_y]$ and $[\kappa_{xx}, \kappa_{yy}]$ are drawn from the same distribution as before. We investigate the effect of different parameter settings than those we used during training to validate the generalizability of the neural scheme. To study the effect of different $\theta$, we use the original test set. We generate two additional test cases varying one parameter at a time: a) $\Delta t=0.12$, and b) $\Delta x =0.049$. We implement all differential operators using the convolutional kernel. The convolutions have kernel size 3 $\times$ 3 $\times$ 3 and are unique for each differential operator. Following \cite{hsieh2018learning}, we use a three-layer convolutional neural network to model each of the $H_i$ with constant padding of one in between. The network is trained using Adam Optimizer with a learning rate $1\times10^{-3}$, betas 0.9 and 0.99 for 20 epochs. The total training time is $\sim$6 hours in a Quadro P6000 GPU.

\begin{figure}[t!]
\centering
\includegraphics[width=0.32\textwidth]{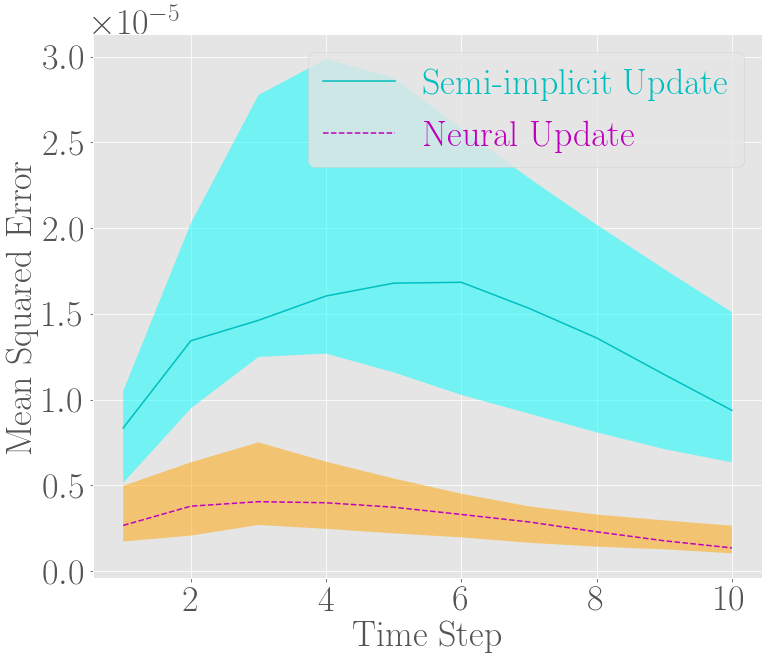}
\includegraphics[width=0.32\textwidth]{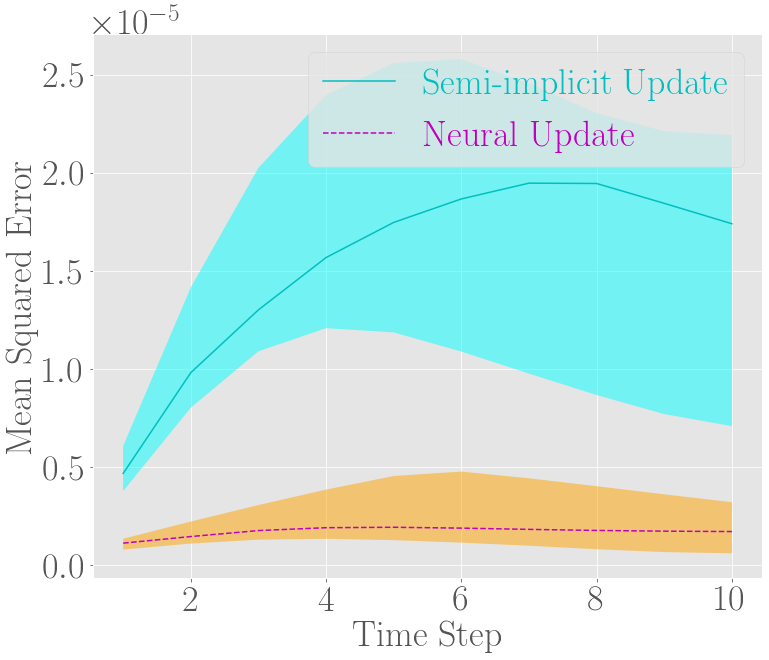}
\includegraphics[width=0.32\textwidth]{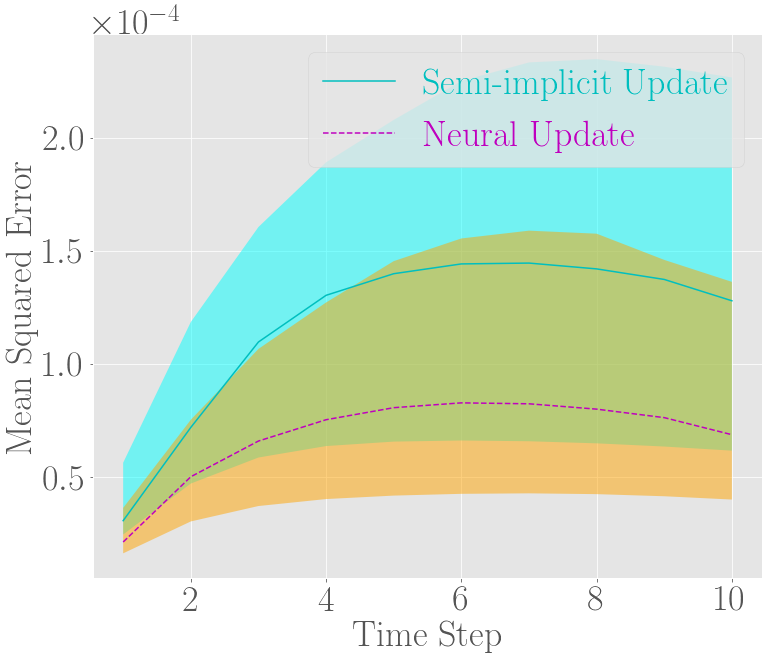} \\
(a) $\theta=0.75$ \hspace{8em} (b) $\Delta t= 0.12$ \hspace{7em} (c) $\Delta x= 0.049$
\caption{ (a), (b), and (c) shows the mean-squared error (between FEniCS solution and semi-implicit scheme and neural scheme) vs a fixed number of time steps plot for different $\theta, \Delta t$, and $\Delta x$ parameters during test time, respectively. The banded curves indicate the 25\% and 75\% percentile of the normalized errors among 20 test samples.}
\label{fig:fig1_exp1}
\end{figure}

\subsubsection{Results}

From Fig.~\ref{fig:fig1_exp1}, we see that error from the neural scheme is less compared to the error from the semi-implicit solution for all three different test sets, with varying $\theta$, $\Delta t$, and $\Delta x$ respectively. Note that we sample the value of $[v_x, v_y] \mbox{ and } [\kappa_{xx}, \kappa_{yy}]$ from the same ranges as the training set, however the exact values are different. We observe that one neural solver iteration takes approximately twice the time of one semi-implicit iteration. Hence, for a fair comparison, we use 10 iterations per time step for our neural solver compared to 25 iterations for the semi-implicit solver. Thus, this experiment affirms our hypothesis that the neural solver is more accurate compared to the semi-implicit solution while keeping generalizability to other PDE hyper-parameter settings at the same time.
\begin{figure}[ht!]
	\centering
	\footnotesize
	\includegraphics[width=0.48\textwidth]{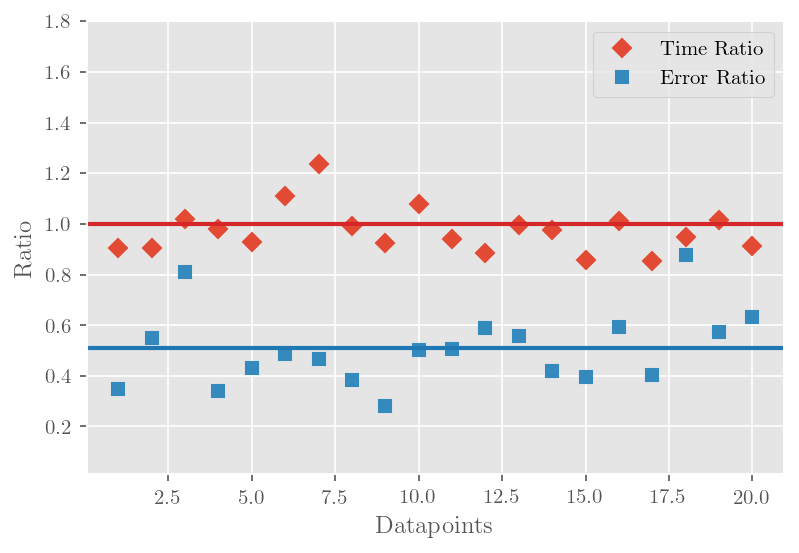}
	\includegraphics[width=0.48\textwidth]{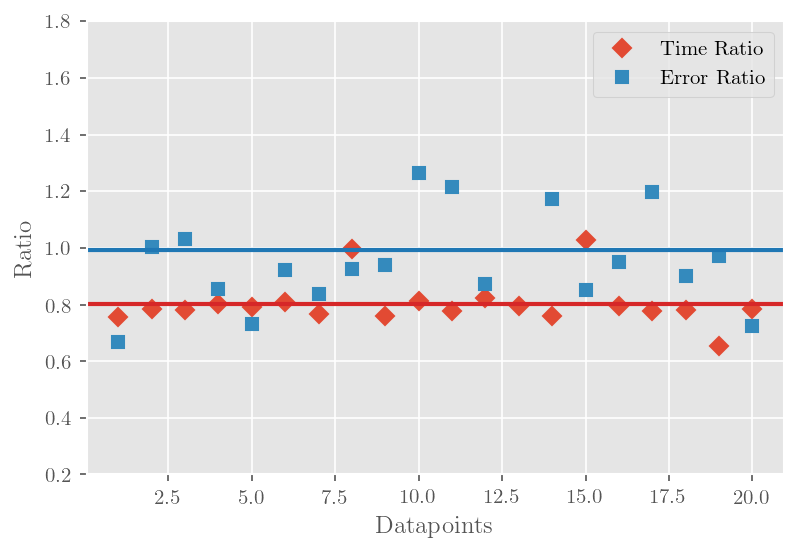}
	(a) Performance gain (blue$<$1) for same \hspace{2cm} (b) Runtime gain (red$<$1) for same\\
	computation budget (red$\sim$1) \hspace{4cm} error tolerance (blue$\sim$1)
	\caption{Performance comparison for the 2D experiment: We plot the runtime and mean-squared error ratio (y-axis) between the neural method and the semi-implicit method for the test samples (x-axis). On the left, we make both neural and semi-implicit solvers run for the same time and compare their errors. We observe that for a given runtime, the neural solver reduces the error by 49.3\%. On the right, we set the same mean-squared error and compare their runtime. We observe that for a given error acceptance, the neural solver is 19.2\% faster.}
	\label{fig:ratios_2d}
\end{figure}

Fig. \ref{fig:fig2_exp1} shows a typical test sample solution from the neural and semi-implicit schemes against our reference FEniCS solution. Qualitatively, we find that the neural solver is producing an accurate solution with a lesser number of iterations, which suggests the learned CNN achieves faster convergence.

For runtime comparison, we chose CPU time because we also want to compare against FEniCS, which runs on CPU. We compare the run time for the neural solver (10 iterations per time step) and semi-implicit scheme (25 iterations per time step). Note that here we only compare the test time of the neural solver after it is trained. The experiments are conducted on an Intel Xeon W-2123 CPU @ 3.60GHz, with code running on one of the four cores. Fig. \ref{fig:ratios_2d} shows the ratio of runtime and error between neural solver and semi-implicit solver. We performed two experiments: we make both neural and semi-implicit solvers run for the same time and compare their error; we set the same average error and compare the runtime of the solvers. From Fig. \ref{fig:ratios_2d}a, we can clearly see that given a fixed runtime neural solver produces 49.3\% error reduction, whereas for a given error acceptance the neural solver is 19.2\% faster, Fig. \ref{fig:ratios_2d}b. We also observe that the trained neural solver takes circa 0.0148s compared to 0.0141s for the semi-implicit scheme, whereas the FEniCS solution takes 3.19s for machine precision convergence.

%%%%%%%%%%%%%%%%%%%%%%%%%%%%%%%%%%%%%%%%%%%%%%%%%%%%%%%%%%%%%%%%%%%%%%%%%%%%%%%%%%%%%%%
\subsection{3D Non-linear PDE with Neumann Boundary}
\label{exp:3d}
We consider a Fisher–Kolmogorov equation for modelling tumor cell density $u$ in human brains as used, e.g., in \cite{lipkova2019personalized}. This equation consists of a reaction and a diffusion term, which is as follows
\begin{equation}
\frac{\partial \mathscr{u}}{\partial t} = \nabla \cdot (\kappa\nabla \mathscr{u}) + \rho \mathscr{u}(1-\mathscr{u});
\quad \textrm{subject to} \quad \nabla \mathscr{u} \cdot \eta = 0 \quad \textrm{at} \quad \Gamma
\end{equation}
where $\kappa$ is given by:
\begin{equation} 
\kappa_i = 
\begin{cases}
p_{w_i} \kappa_w + p_{g_i} \kappa_g & \text{for \textit{i} in } \Omega \\
0 & \text{elsewhere}
\end{cases}
\end{equation}
$p_{w_i}$ and $p_{g_i}$ denote the \% of the white and gray matter tissue at the $\textit{i}^{th}$ voxel, respectively. The constants $\kappa_w$ and $\kappa_g$ describe tumor infiltration rate in white and gray matter, respectively and it is assumed that $\kappa_w = 10\times \kappa_g$. The unit for $\kappa$ is [$\frac{cm^2}{day}$]. $\rho$ parameterizes the proliferation rate (the number of cells that divide per day). Its unit is $\frac{1}{day}$.

It can already be seen that this equation differs in several ways from the previously considered 2D equation and other commonly used equations in the existing literature that aim at leveraging neural networks in their solver strategy. The Fisher–Kolmogorov equation is a non-linear time-dependent PDE with spatially dependent PDE parameters subject to Neumann boundary condition, which to the best of our knowledge, is being used in this paper for the first time when leveraging a neural solver.  Also, it is important to point out that the human brain possesses very irregular geometry, which in turn offers a wide variety of simulation domains. Although there is no theoretical guarantee that the learned neural solver will generalize to an unseen domain, e.g., arising from the brain anatomy of a new patient, we experimentally show that it does.

\begin{figure}[t!]
\centering
\includegraphics[width=0.97\textwidth]{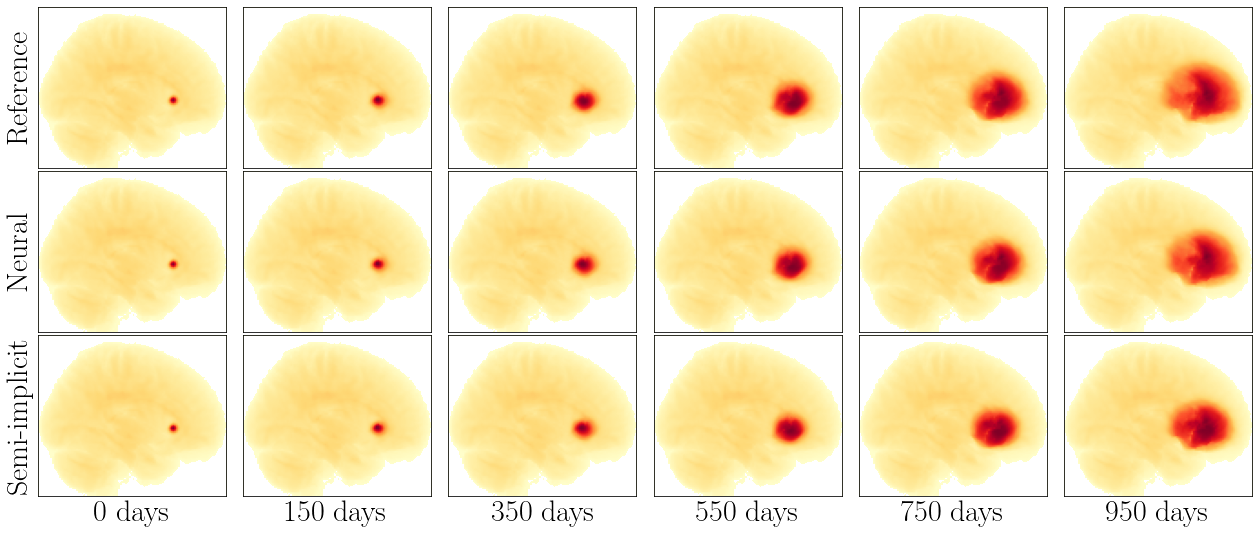}
\includegraphics[width=0.97\textwidth]{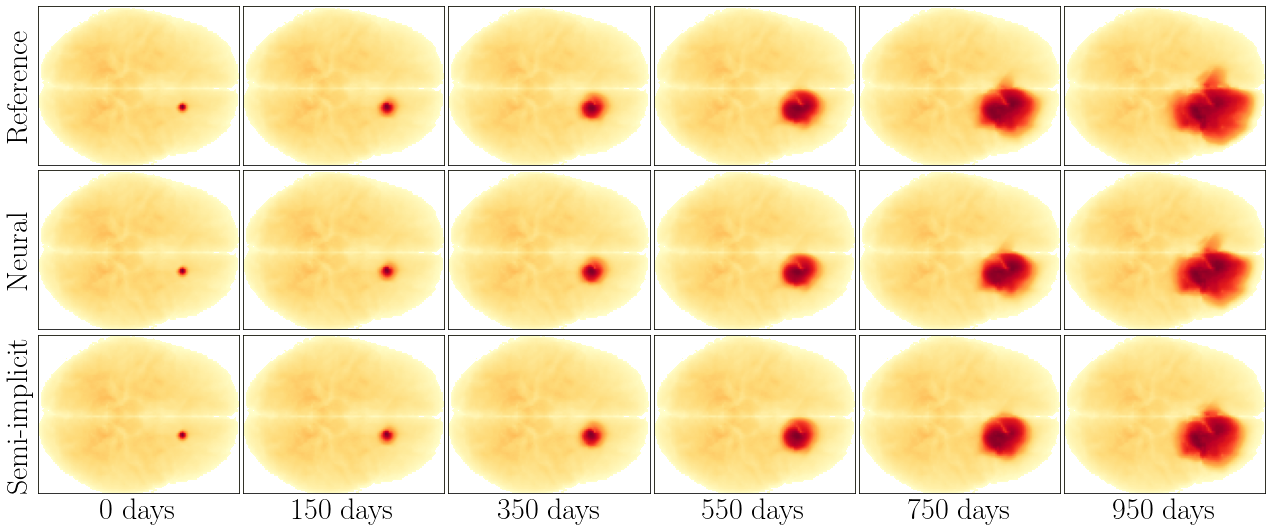}
	\caption{2D mean intensity projection of the 3D tumor cell density (in red-orange) $u$ (c.f. Eq 18) overlayed with brain tissue map (in yellow) for the neural scheme and a semi-implicit scheme against the reference solution for two test sequences over times. All methods use the same initial- and boundary conditions. The neural update shows consistently faster convergence than the semi-implicit one.}
	\label{fig:qualitative}
\end{figure}

Following the diffuse domain approach in Sec. \ref{sec:phase_field}, we reformulate the problem as
\begin{align}
\phi_i &= p_{w_i} + p_{g_i} \\
\frac{\partial \phi \mathscr{u}}{\partial t} &= \phi\nabla \cdot ( \kappa\nabla \mathscr{u}) + \nabla \phi \cdot \nabla \mathscr{u}+ \phi \rho \mathscr{u}(1-\mathscr{u}) 
\end{align}
This allows us to recast the Neumann boundary condition to a Dirichlet boundary condition given by:
\begin{equation}
\mathscr{u} = 0 \quad \text{at} \quad \Gamma_1
\end{equation}
where $\Omega_1$ is the entire $129 \times 129 \times 129$ voxels volume. To adapt this reformulation to fit into our update rule, we write,
\begin{equation}
	\phi\frac{u_{t+\Delta t} - u_t}{\Delta t} = \theta (\phi\nabla \cdot ( \kappa\nabla u_{t + \Delta t}) + \nabla \phi \cdot \nabla u_{t + \Delta t})) + (1-\theta) (\phi\nabla \cdot ( \kappa\nabla u_{t}) + \nabla \phi \cdot \nabla u_{t})) + \phi r(u_{t})
\end{equation}

\begin{figure}[t!]
	\centering
	\includegraphics[width=\textwidth]{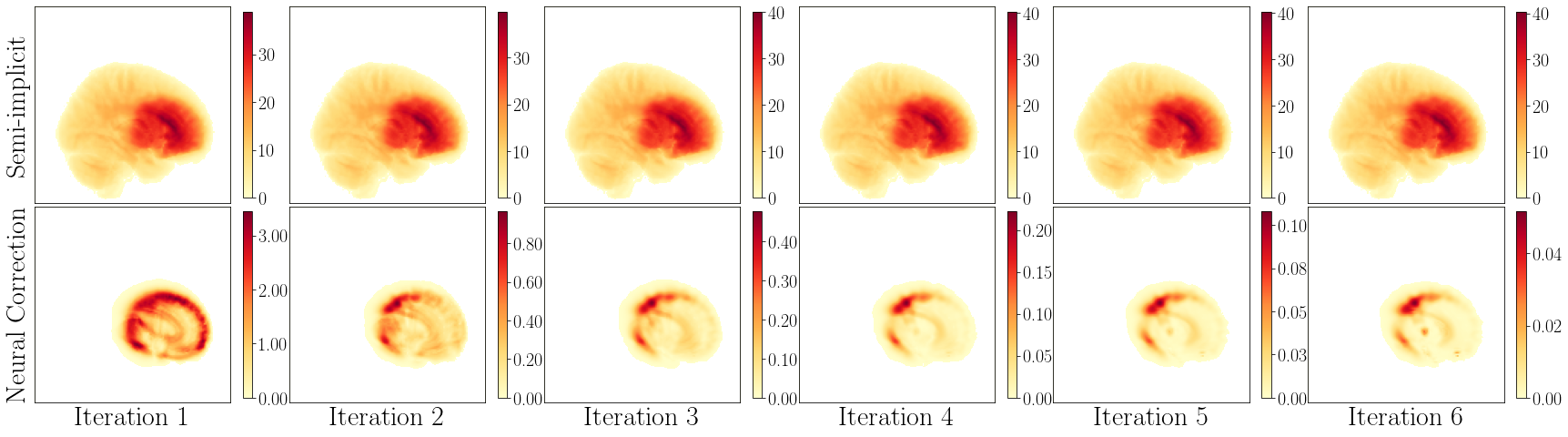}
	\caption{2D mean intensity projection of an example of the two terms constituting a PDE solution (Eq. 11) plotted from left to right for the first 6 iterations of a single time step. On the top, we see the domain (the brain in yellow overlay) and the semi-implicit part of the tumor density (in red-orange), $\Psi(u)$. At the bottom, we see the neural correction component of the tumor density (without brain overlay), $G \left( \sum_{i=1}^{N} \Lambda_i H_i w \right)$.  It shows that the neural correction makes a significant impact in the initial iterations, which helps the learned solver converge faster.}
	\label{fig:correction}
\end{figure}

\subsubsection{Data Generation: }
\label{sec:datagen}
For our training data, we randomly generated 72 initial conditions $u_0$ with different $\rho$ and $\kappa_w$ at random initial tumor locations. To acquire the reference data, we let a semi-implicit solver run to convergence for each time step. Here the reference solution means, in this case, the fixed point $u^*$ of the iterator on which the neural solver is built (in this case, the semi-implicit solver). We obtain the reference solution by running the semi-implicit solver in implicit mode ($\theta=1$ for maximum convergence speed) until convergence. This way, the ground truth data need not necessarily be simulated and stored beforehand and can be computed on the fly during training. Per initial condition, we run the training for 19 time-steps with $\Delta t=50$ days, $\Delta x=0.15, \Delta y=0.18, \Delta z=0.15$ mm. Thus, in the end, we have 20 time points for each sample.

\begin{figure}[ht!]
	\centering
	\footnotesize
	\includegraphics[width=0.48\textwidth]{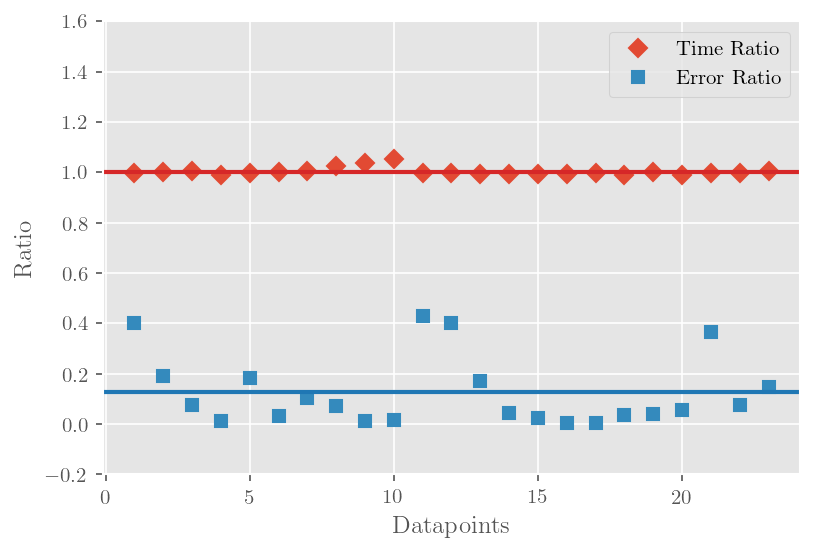}
	\includegraphics[width=0.48\textwidth]{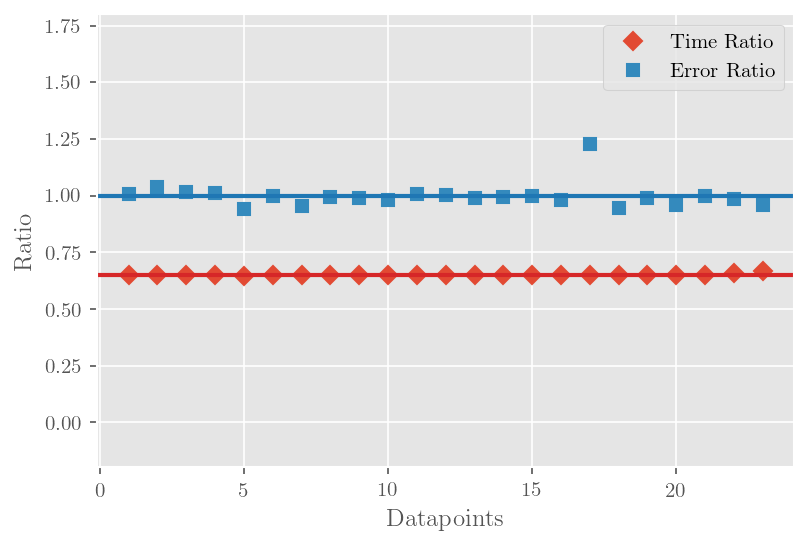}
	\\
	(a) Performance gain (blue$<$1) for same \hspace{2cm} (b) Runtime gain (red$<$1) for same\\
	computation budget (red$\sim$1) \hspace{4cm} error tolerance (blue$\sim$1)
	\caption{Performance comparison for the 3D experiment: We plot the runtime and mean-squared error ratio (y-axis) between the neural method and the semi-implicit method for the test samples (x-axis). On the left, we make both neural and semi-implicit solvers run for the same time and compare their errors. We observe that for a given runtime, the neural solver reduces the error by 87\%. On the right, we set the same mean-squared error and compare their runtime. We observe that for a given error acceptance, the neural solver is 35\% faster.}
	\label{fig:ratios_3d}
\end{figure}

\subsubsection{Experimental Details}

We assign randomly selected 49 samples as training and 23 samples as the test set. After we generated our training data and found our reference, we run the neural solver on a random number of iterations (between 5 and 10) per time step. The cost is evaluated as the mean squared error between the previously calculated fixed point (the reference solution) and the neural output, as given in Equation \eqref{eq:costfunc}. We implement all differential operators using the convolutional kernel. The convolutions have kernel size 3 $\times$ 3 $\times$ 3 and are unique for each differential operator. The network is trained using Adam Optimizer with a learning rate of 1e-3, betas 0.9, and 0.99.  The training error converged after 13 epochs.

\subsubsection{Results}
\label{sec:results}

To better understand the mechanism of the solver with visual illustration, we look at few iterations of an example solution of the PDE, Fig. \ref{fig:correction}. Importantly, Fig. \ref{fig:correction} highlights that the neural part works as a corrective rather than a prediction of its own. The neural net has learned to correct the errors from the semi-implicit iterator, so that the fixed point is reached faster. Note, however, that the speed of growth naturally depends significantly on the parameters. Also, an important observation here is that the correction terms are stronger in the initial iterations since the discrepancy between the fixed point and the current solution is larger in the initial iterations.

Fig. \ref{fig:qualitative} shows a qualitative visual comparison between the neural solver and the semi-implicit solver in contrast to the reference solution. While both the neural solver and the semi-implicit solver took the same runtime, the neural solver shows faster convergence as the time step increases.

Figure \ref{fig:ratios_3d} shows the ratios between the semi-implicit solver and the neural solver for their respective computation time and the mean absolute errors for the test data. As in the 2D case, we performed two experiments testing the solvers' error upon a fixed running time and the solver's running time upon a fixed error level. We observe that for a given runtime, the neural solver reduces the error by 87\%, Fig. \ref{fig:ratios_3d}a. In addition, for a given error acceptance the neural solver is 35\% faster, Fig. \ref{fig:ratios_3d}b.

In our experiments, the test set consists of parameters and computation domain, which were not seen during the training time. This confirms that the solver is able to generalize over other parameters and computation domains.

\subsection{Discussion}
\label{exp:disc}
From two experiments, we draw the following observations.
\begin{itemize}
    \item Both experiments confirm the theoretical convergence property of the neural solver in different PDE parameter settings.
    \item We observe that our solver is well applicable to the given advection-diffusion and reaction-diffusion equations in 2D and 3D, which indicates its generalizability to other linear PDEs as well.
    \item In the 3D experiment, we have very different computation domains in a different region of brain tissue. Hence, we empirically show that our solver generalizes to the arbitrary shape of the boundary.
    \item The reaction-diffusion experiment shows that our model can handle the specific class of non-stiff non-linear PDEs really well.
    \item We observe significant speedup compared to the traditional solver in both experiments. We observe this to be more pronounced in the 3D setting than in 2D. We hypothesize that in 3D, the neural solver can leverage the domain-specific spatially varying PDE parameters more efficiently in the learned PDE solutions.
\end{itemize}

%% file: 04_appendices.tex
\section{Proofs}
\label{appendixA}
\thmone*
\begin{proof}
The spectral radius of $L$ can be bounded by: 
\begin{align}
\rho(L) \leq ||L|| &= \left| \left| G \sum_{i = i}^{N} \Lambda_i (D_i - d_i I)\right| \right| \nonumber\\
 &\leq ||G|| \sum_{i = i}^{N} || \Lambda_i || \text{ } || (D_i - d_i I)|| \nonumber\\
 &= \sum_{i = i}^{N} || \Lambda_i || \text{ } || (D_i - d_i I)||; \text{ [since $||G|| = 1$]} \nonumber
\end{align}

Thus given $||\Lambda_i|| < \frac{1}{\sum_{j = i}^{N} || (D_j - d_j I)||}, \forall i \leq N$ we have $\rho(L) < 1$. We can therefore conclude that the newly introduced iterator from Eq. \ref{eq:new_it} is a valid \textit{fixed-point iterator} for a proper choice of $\Delta x, \Delta t$ and $\theta$.
\end{proof}

\thmthree*
\begin{proof}
The proof is similar to that of \eqref{thm3}. The spectral norm $|| . ||$ is convex form the subadditive property and $L'$ is linear in $\mathbb{H}$. To prove that it is open, observe that $||.||$ is a continous function, so $(L + G \sum_{i = 1}^N \Lambda_i H_i (L - I))$ is continous in $\mathbb{H}$. Given $\rho(L')$, the set of $\mathbb{H}$ is the preimage under this continous function of $(0, 1 - \zeta)$ for some $\zeta > 0$, and the inverse image of the open set $(0, 1 - \zeta)$ must be open. 
\end{proof}

\lemmafour*
\begin{proof}
Considering the spectral norm of $L'$ and invoking product and triangular inequality of norms, we obtain the following tight bound: \begin{align}
	\rho(L') \leq ||L'|| &= \left| \left|G \sum_{i=1}^{N} \Lambda_i (D_i - d_i I - H_i) + G \sum_{i=1}^{N}(\Lambda_i H_i)L \right| \right| \nonumber\\
	&\leq ||G|| \sum_{i=1}^{N} ||\Lambda_i|| \text{ } ||D_i - d_i I - H_i|| + ||G|| \sum_{i=1}^{N}||\Lambda_i|| ||H_i|| ||L|| \nonumber\\
	&< \sum_{i=1}^{N} ||\Lambda_i|| (||D_i - d_i I - H_i|| + ||H_i||) \quad [||G|| = 1, ||L|| < 1]\nonumber
\end{align} Given $||L|| < 1$ we have $||\Lambda_i|| < \frac{1}{\sum_{j = i}^{N} || (D_j - d_j I)||}, \forall i \leq N$, hence  
\begin{align}
	\rho(L') &< \frac{1}{\sum_{j = i}^{N} || (D_j - d_j I)||} \sum_{i = 1}^N(|| D_i - d_i I - H_i|| + ||H_i||)  \nonumber\\
	&<1;~~[
\mbox{Invoking triangular inequality}]\nonumber
\end{align}
\end{proof}

\corofive*
\begin{proof}
From Theorem \eqref{thm1} and Lemma \eqref{lemma2} we know that our iterator is valid if and only if $\rho(L') < 1$. From Lemma \eqref{lemma4} the upper bound of the spectral norm of the iterator depends only on $\mathbb{D}$ and $\mathbb{H}$ given $||\Lambda_i|| < \frac{1}{\sum_{j = 1}^{N} || (D_j - d_j I)||}, \forall i \leq N$ Nonetheless, for any matrix the spectral radius is upper bounded by its spectral norm. Thus, if the iterator is valid for some $u_t', b_{t+\Delta t}', \mathbb{A}, \Delta x', \Delta t',$ and $\theta'$ then it is feasible for any choice of $u_t, b_{t+\Delta t}, \mathbb{A}, \Delta x, \Delta t,$ and $\theta$ that satisfy the constraints.
\end{proof}